\documentclass[a4paper,ejs,preprint,reqno]{imsart}

\usepackage{amsthm,amsmath}
\usepackage[authoryear]{natbib}

\usepackage[psamsfonts]{amssymb}

\RequirePackage[colorlinks,citecolor=blue,urlcolor=blue,linktocpage=true]{hyperref}

\usepackage{dsfont}

\DeclareMathOperator{\sign}{sgn}

\def\hat{\widehat}
\newcommand{\RR}{\mathbb{R}}
\newcommand{\bX}{\boldsymbol X}
\newcommand{\bfE}{\mathbf E}
\newcommand{\calX}{\mathcal X}
\newcommand{\calY}{\mathcal Y}
\newcommand{\calH}{\mathcal H}

\newcommand{\1}{\mathds 1}

\begin{document}

\begin{frontmatter}

\title{Discussion on the paper: Hypotheses testing by convex optimization by Goldenshluger, Juditsky and Nemirovski}
\runtitle{{Discussion on ``Hypotheses testing by convex optimization''}}

\begin{aug}
\author{\fnms{Arnak S.} \snm{Dalalyan},\ead[label=e1]{arnak.dalalyan@ensae.fr}
\ead[label=u1,url]{http://arnak-dalalyan.fr/}}

\runauthor{Dalalyan, A.S.}

\affiliation{ENSAE ParisTech-CREST}

{\renewcommand{\addtocontents}[2]{}
\address{{3, Avenue Pierre Larousse,}\\
{92240 Malakoff,}
{France} \\
\printead{e1}\\
\printead{u1}}
}
\end{aug}


%
%


\end{frontmatter}

\maketitle


This is an exciting piece of work. I agree with the authors that developing computationally tractable
methods for hypotheses testing is an important problem in statistics that have received little attention to date.
In this discussion, I would like to put the emphasis on three points presented in the paper
under discussion that are of particular interest.

\subsection*{Connection with the statistical learning theory}

The idea of convexification of the loss function in order to construct computationally tractable procedures
has been widely used in statistical learning theory \citep{Zhang1}. In this part of the discussion, I would like
to share some thoughts about the similarities of the two approaches.

To this end, let me briefly recall the principle of loss convexification in the problem of binary classification.
One observes $n$ iid pairs $\{(X_i,Y_i)\}_{i=1,\ldots,n}$ drawn from an unknown distribution $P$ on the
product space $\calX\times\calY$ with $\calY=\{-1,+1\}$ and the goal is to design a prediction rule
$g:\calX\to\calY$ with the smallest possible misclassification error rate
\begin{align}\label{risk}
R_P(g) = \bfE_P[\1(Y\not = g(X))]   = \bfE_P[\1(-Yg(X)\ge 0)].
\end{align}
The convexification is achieved in two steps. First, the classification risk is replaced by the $\phi$-risk
\begin{align}\label{riskA}
A_P(g) = \bfE_P[\phi(-Yg(X))],
\end{align}
where $\phi:\RR\to\RR$ is a convex function often referred to as the convex surrogate loss. Second, the set of ``pure''
classification rules $g:\calX\to\calY$ is extended to ``generalized'' rules $h:\calX\to\RR$ with the convention
that the predictions furnished by $h$ and $\sign(h)$ are the same. The $\phi$-risk is accordingly extended
to all generalized prediction rules: $A_P(h)=\bfE_P[\phi(-Yh(X))]$. As a consequence of this construction,
if $\calH$ is a convex subset of the set of all measurable functions from $\calX$ to $\RR$, then the computation of the empirical
risk minimizer (ERM)
\begin{align}\label{ERM}
\hat h_{n,\calH} \in\text{arg}\min_{h\in\calH} \frac1n\sum_{i=1}^n \phi(-Y_ih(X_i))
\end{align}
amounts to solving a convex program. The most common choices for the function $\phi$ are the hinge loss
$\phi(u) = (1+u)_+$, the exponential loss $\phi(u)= e^u$ and the logistic loss $\phi(u) = \log(1+e^u)$.

Let me turn now to the problem of testing two hypotheses  $\Theta_{0}$ and $\Theta_{1}$ based on $n$ observations
$X_1,\ldots,X_n$ independently drawn from a distribution $P_{\theta^*}$ on $\calX$.  Let $\Theta = \Theta_0\cup \Theta_1$ and
$s:\Theta\to\{\pm1\}$ be the function that equals $-1$ on the set $\Theta_0$ and $+1$ on the set $\Theta_1$.
The usual loss of a pure test $T:\calX^n\to\{\pm 1\}$ associated with a sample $\bX^{(n)}=(X_1,\ldots,X_n)$ drawn from
$P_{\theta^*}^n:=P_{\theta^*}^{\otimes n}$ is
$$
\ell\big(T(\bX^{(n)}),\theta^*\big) = \1\big(T(\bX^{(n)})\not = s(\theta^*) \big) = \1(-T(\bX^{(n)})s(\theta^*)\ge 0).
$$
The corresponding risk is $R_{P_{\theta^*}}(T) = \bfE_{{\theta^*}}\big[\ell\big(T(\bX^{(n)}),\theta^*\big)\big]$
and the worst case risk is
\begin{align}\label{eps}
\epsilon(T) &= \sup_{\theta\in\Theta}R_{P_{\theta^*}}(T) = \sup_{\theta\in\Theta} P_{\theta}\big(T(\bX^{(n)})\not = s(\theta)\big)\nonumber\\
&=\sup_{\theta\in\Theta} \bfE_{{\theta}}[\1(-T(\bX^{(n)})s(\theta)\ge 0)].
\end{align}
Comparing (\ref{eps}) with (\ref{risk}), one can see some clear similarities between the problems of finding binary predictors $g$
minimizing the misclassification error rate and that of finding testing procedures $T$ minimizing the worst case error rate $\epsilon(T)$.
In both problems the decision rules form a nonconvex set and the performance measure is defined as the expected loss for a nonconvex
loss function  (the Heaviside step function). However, there is an important difference consisting in the fact that---contrary to
(\ref{risk})---the expectation at the right-hand side of (\ref{eps}) does not admit an empirical counterpart that is easily computable
from the sample.  Therefore, even if one applies the aforementioned two steps of convexification, this does not readily yield a test procedure
computable by solving a convex program (in the spirit of (\ref{ERM})).

Elaborating on these ideas,
one can define the following  convexified strategy for testing the hypothesis $\Theta_0$ against $\Theta_1$. Given a convex subset $\calH$ of the set of measurable functions
from $\calX^{n}$ to $\RR$ and a convex loss $\phi:\RR\to\RR$, define
\begin{align}\label{test}
\hat h_{n,\calH}^{\phi} \in \text{arg}\min_{h\in\calH} \sup_{\theta\in\Theta} G_\phi(h,\theta),
\qquad
G_\phi(h,\theta)=
\bfE_{\theta}\big[\phi\big(-h(\bX^{(n)})s(\theta)\big)\big].
\end{align}
In this ``saddle-point'' formulation, the outer minimisation problem has the attractive property of being convex: it has a convex feasible set and a convex cost function.
Unfortunately, in general, the inner maximization problem is not concave and there is no particular reason to expect that it can be efficiently solved
for any given $h$ when the dimensionality of $\theta$ is large. To circumvent this drawback, the authors 
had the ingenious idea to combine the following three facts:
\begin{itemize}
\item the saddle point of $G(h,\theta)$ coincides with the saddle point of $\log G(h,\theta)$,
\item when the model $\{P_\theta:\theta\in\Theta\}$  belongs to an exponential family, it is natural to
choose $\calH$ as the linear span of the sufficient statistics: $\calH_0 = \text{Span}(S_j:j=1,\ldots,m)$,
\item for some statistical models\footnote{It could be helpful to mention that the concavity property holds for the usual 
parameterization and does not hold for the parameterization in terms of the natural parameters in the sense of exponential families.} 
belonging to an exponential family, for every $h\in\calH_0$,
the mapping $\theta\mapsto \log\big(\bfE_{\theta}[ \exp(-h(\bX^{(n)}))]\big)$  is concave.
\end{itemize}
This leads to the test procedure
\begin{align}\label{T-exp}
\hat h_{n,\calH_0}^{\rm exp} \in \text{arg}\min_{h\in\calH_0}\ \sup_{\theta\in\Theta}\ \ \log G_{\rm exp}(h,\theta),
\qquad G_{\rm exp}(h,\theta)=\bfE_{\theta}\big[e^{-h(\bX^{(n)})s(\theta)}\big].
\end{align}
The final step of construction aims at convexifying the feasible set of the inner maximization problem. In the case when $\Theta = \Theta_0\cup \Theta_1$
with convex sets $\Theta_0$ and $\Theta_1$, this aim is achieved by replacing $\sup_{\theta\in\Theta}\ \log G_{\rm exp}(h,\theta)$ by the expression
$\sup_{(\theta,\bar\theta)\in\Theta_0\times\Theta_1}\  \log G_{\rm exp}(h,\theta) + \log G_{\rm exp}(h,\bar\theta)$, which does not impact the error of testing too
much in view of the inequalities
\begin{align*}
\sup_{\theta\in\Theta}\ \log G_{\rm exp}(h,\theta) &\le \sup_{(\theta,\bar\theta)\in\Theta_0\times\Theta_1}\  \big\{\log G_{\rm exp}(h,\theta) +
\log G_{\rm exp}(h,\bar\theta)\big\} \\
&\le 2\sup_{\theta\in\Theta}\ \log G_{\rm exp}(h,\theta).
\end{align*}
An important remark to be made here is that---in the case of exponential loss $\phi$---taking the logarithm of $G_\phi$
does not break the convexity with respect to $h$. So, in this notation, the test proposed and studied by the authors is
\begin{align}\label{T-exp1}
\tilde h_{n,\calH_0}^{\rm exp} \in \text{arg}\min_{h\in\calH_0}\ \sup_{(\theta,\bar\theta)\in\Theta_0\times\Theta_1}\
\big\{\log G_{\rm exp}(h,\theta) + \log G_{\rm exp}(h,\bar\theta)\big\}.
\end{align}
I believe that these explanations shed some additional light on the construction proposed in Theorem~2.1 of the paper under discussion.
This also raises several questions that might be interesting to investigate in the future. In particular, a compelling question is to
characterize the set of surrogate loss functions $\phi$ that lead to computationally tractable testing procedures and for which the testing
error rate remains small. Another question is the possibility to deal with test (\ref{T-exp}) directly, without using the final step of
convexification. At a heuristic level, the risk of $\hat h_{n,\calH_0}^{\rm exp}$ should be smaller than that of $\tilde h_{n,\calH_0}^{\rm exp}$.
Therefore, the advantage of the latter would be only computational tractability. I wonder if it is possible to efficiently compute
the test $\hat h_{n,\calH_0}^{\rm exp}$, despite the lack of convex-concavity of the cost function, exploiting the facts that (a)
for every $h$, the sup of $\log G_{\rm exp}(h,\theta)$ over $\Theta$ can be efficiently computed, and (b) for every $\theta$, the minimum
of $\log G_{\rm exp}(h,\theta)$ over $\calH_0$ can be efficiently computed as well.

\subsection*{Reduction to testing simple hypotheses}

The definition of the test given by the authors in Theorem 2.1, see also Eq.\ (\ref{T-exp1}) above, is well suited for the computational
purposes but, in my opinion, has the inconvenience of hiding the main reason why the proposed test is a natural one to use in the setting
under consideration. In fact, the proposed test can be alternatively defined as follows: in order to distinguish between two
(convex) hypotheses $\Theta_0$ and $\Theta_1$ based on a sample $\bX\sim P_{\theta^*}$,
\begin{enumerate}
\item Determine the two closest points $\theta_0\in\Theta_0$ and $\theta_1\in\Theta_1$ in terms of the Hellinger distance between the corresponding
distributions (in other terms, find the two representers $P_{\theta_0}$ and $P_{\theta_1}$ in the families $\{P_\theta:\theta\in\Theta_0\}$ and
$\{P_\theta:\theta\in\Theta_1\}$ that are the hardest to distinguish). This step is completely data independent.
\item Apply the standard likelihood-ratio test to the problem of choosing among two simple hypotheses $H_0: \theta = \theta_0$ versus
$H_1:\theta=\theta_1$.
\end{enumerate}
The equivalence of these two definitions follows from the proof of Theorem 2.1, see Eq.\ (52). In Section 2.3.2, this interpretation is presented
for the discrete observation scheme. At a conceptual level, it is important to underline that the same interpretation holds true in the general case
as well. However, from a practical point of view, the definition given in the paper is more convenient than the foregoing one since the first step of
the latter, generally, is not computationally tractable.

\subsection*{Testing error for inexact solutions}

As it is judiciously noted by the authors, in many practical situations, the exact computation of the saddle point in (\ref{T-exp1}) can 
not be performed. Then, one relies on an approximation of the saddle point and it is a central task to assess how this approximation error impacts
the error of testing. I find it relevant to measure the error of approximation in terms of the magnitude of violation of first-order optimality
conditions (see, for instance, Eq.\ (8) of the paper under discussion). In such a context, the authors establish upper bounds on the error
of the test based on an approximate solution to the saddle point problem. For example, in the case of the Gaussian observation scheme explored
in Section 2.3.1, it is shown that the worst-case error rate of the test based on the exact solution is
\begin{equation}\label{bound0}
\varepsilon_* = 1-\Phi\Big(\frac12\big\|\Sigma^{-1/2}(\theta_0-\theta_1)\big\|_2\Big),
\end{equation}
where $\Phi$ is the cumulative distribution function of the standard normal distribution and $(\theta_0,\theta_1)$ is the second argument
of the solution to the saddle point problem. On the other hand, when an inexact solution $(\tilde\theta_0,\tilde\theta_1)$ is used, with
an approximation error bounded by $\delta>0$, the worst-case error rate satisfies (see Eq. (9)):
$$
\tilde \varepsilon \le 1-\Phi\Big(\frac12\big\|\Sigma^{-1/2}(\tilde\theta_0-\tilde\theta_1)\big\|_2-\frac{\delta}{\|\Sigma^{-1/2}(\tilde\theta_0-\tilde\theta_1)\|_2}\Big).
$$
In my opinion, it is worth complementing this upper bound by another one that involves only the exact solution $(\theta_0,\theta_1)$ and, therefore,
makes it easier to compare the two errors $\varepsilon_*$ and $\tilde\varepsilon$. In the case of Gaussian observation scheme, this can be easily done.
In fact, one can deduce from the first-order exact and approximate optimality conditions that
\begin{equation}\label{bound1}
\|\Sigma^{-1/2}(\theta_0-\theta_1)\|_2-\sqrt{\delta}\le \|\Sigma^{-1/2}(\tilde\theta_0-\tilde\theta_1)\|_2\le \|\Sigma^{-1/2}(\theta_0-\theta_1)\|_2+\sqrt{\delta}
\end{equation}
Since the Gaussian cdf is increasing, we infer from this inequality that
$$
\tilde \varepsilon \le 1-\Phi\Big(\frac12\big\|\Sigma^{-1/2}(\theta_0-\theta_1)\big\|_2-\frac{\sqrt{\delta}}{2}-
\frac{\delta}{\|\Sigma^{-1/2}(\theta_0-\theta_1)\|_2-\sqrt{\delta}}\Big).
$$
An even more elegant bound can be obtained if the normalized approximate optimality condition is used:
$\forall (\theta,\bar\theta)\in \Theta_0\times\Theta_1$, it holds
$$
(\tilde\theta_1-\tilde\theta_0)\Sigma^{-1}(\theta-\tilde\theta_0) +
(\tilde\theta_0-\tilde\theta_1)\Sigma^{-1}(\bar\theta-\tilde\theta_1) \le \delta\|\Sigma^{-1/2}(\tilde\theta_0-\tilde\theta_1)\|_2^2.
$$
In this case, inequalities (\ref{bound1}) take the form
\begin{equation}\label{bound1'}
\frac{\|\Sigma^{-1/2}(\theta_0-\theta_1)\|_2}{1+\sqrt{\delta}}\le \|\Sigma^{-1/2}(\tilde\theta_0-\tilde\theta_1)\|_2\le
\frac{\|\Sigma^{-1/2}(\theta_0-\theta_1)\|_2}{1-\sqrt{\delta}}
\end{equation}
and we get
$$
\tilde \varepsilon \le 1-\Phi\bigg\{\Big(\frac12-\delta\Big)\frac{\big\|\Sigma^{-1/2}(\theta_0-\theta_1)\big\|_2}{1+\sqrt{\delta}}\bigg\}.
$$
This inequality allows for an easy comparison of $\tilde\varepsilon$  and $\epsilon_*$ in the case of Gaussian observations. In the case
of other observation schemes, deriving this type of upper bounds seems to be more challenging and constitutes an interesting avenue of future
research.

{\renewcommand{\addtocontents}[2]{}
\section*{Acknowledgments}
The research of the author is partially supported by the grant Investissements d'Avenir (ANR-11-IDEX-0003/Labex Ecodec/ANR-11-LABX-0047).
}

\bibliographystyle{plainnat}

{\renewcommand{\addtocontents}[2]{}
\bibliography{Literature1}}

\end{document}